\documentclass[11pt]{article}

\usepackage{amsmath,amssymb,amsfonts,amsthm}

\newcommand{\Hom}{\mbox{Hom}\,}
\newcommand{\Ext}{\mbox{Ext}\,}

\newcommand{\Spec}{\mbox{Spec}\,}

\newcommand{\Ass}{\mbox{Ass}\,}
\newcommand{\Supp}{\mbox{Supp}\,}

\renewcommand{\dim}{\mbox{dim}\,}

\newcommand{\q}{\mbox{q}\,}

\renewcommand{\H}{\mbox{H}}
\newcommand{\V}{\mbox{V}}

\newcommand{\lo}{\longrightarrow}

\newcommand{\fa}{\frak{a}}
\newcommand{\fb}{\frak{b}}
\newcommand{\fm}{\frak{m}}
\newcommand{\fp}{\frak{p}}

\begin{document}
\date{}
\title{\bf Associated primes and cofiniteness of local cohomology modules
}

\author{Mohammad T. Dibaei $^{\it ab}$ and Siamak Yassemi$^{\it
ac}$\\
{\small\it $(a)$ Institute for Studies in Theoretical Physics and
Mathematics, Tehran, Iran}\\
{\small\it $(b)$ Department of Mathematics, Teacher Training
University, Tehran, Iran}\\
{\small\it $(c)$ Department of Mathematics, University of Tehran,
Tehran, Iran}} \maketitle

\begin{abstract}

\noindent Let $\fa$ be an ideal of Noetherian ring $R$ and let $M$
be an $R$-module such that $\Ext^i_R(R/\fa,M)$ is finite
$R$-module for every $i$. If $s$ is the first integer such that
the local cohomology module $\H^s_\fa(M)$ is non $\fa$-cofinite,
then we show that $\Hom_{R}(R/\fa, \H^s_\fa(M))$ is finite.
Specially, the set of associated primes of $\H^s_\fa(M)$ is
finite.

Next assume $(R,\fm)$ is a local Noetherian ring and $M$ is a
finite module. We study the last integer $n$ such that the local
cohomology module $\H^n_\fa(M)$ is not $\fm$-cofinite and show
that $n$ just depends on the support of $M$.

\end{abstract}

\baselineskip=18pt

\vspace{.3in}

\noindent{\bf 1. Introduction}

\vspace{.2in}

 Let $R$ be a Noetherian ring, $\fa$ an ideal, and $M$ an
 $R$-module. An important problem in commutative algebra is
 determining when the set of associated primes of the ith local
 cohomology module $\H^i_\fa(M)$ of $M$ is finite. If $R$ is a
 regular local ring containing a field then $\H^i_\fa(R)$ has only
 finitely many associated primes for all $i\ge 0$, cf. [{\bf HS}]
 (in the case of positive characteristic) and [{\bf L}] (in
 characteristic zero). In [{\bf S}] Singh has given an example of
 a Noetherian non-local ring $R$ and an ideal $\fa$ such that
 $\H^3_\fa(R)$ has infinitely many associated primes. On the other
 hand, Brodmann and Lashgari [{\bf BL}; Theorem 2.2] have shown that the first
 non-finite local cohomology module $\H^i_\fa(M)$ of a finite
 module $M$ with respect to an ideal $\fa$ has only finitely many associated
 primes (see also [{\bf KS}]).

 An $R$-module $M$ is said to be $\fa$-cofinite if $\Supp(M)\subseteq\V(\fa)$ and
 $\Ext^i_R(R/\fa,
 M)$ is finite for all $i\ge 0$. Cofiniteness of modules and local cohomology modules
  have been
 studied by many authors, c.f.  [{\bf DeM}], [{\bf Ha2}], [{\bf M1}], [{\bf M2}],
 and [{\bf M3}].
 In section 2, the main aim is to give a generalization of Brodmann and Lashgari's
 result by showing that the $R$-module $\Hom_R(R/\fa,\H^s_\fa(M))$ is finite,
whenever the modules $\H^i_\fa(M)$ are $\fa$-cofinite for all
$i<s$ and $\Ext^j_R(R/\fa, M)$ is finite $R$-module for every
$j\le s$. This result improves [{\bf BL}; Theorem 2.2].

In section 3, we assume that $(R,\fm)$ is a Noetherian local ring.
In [{\bf Ha1}], Hartshorne defined $\q(\fa, R)$ as the least
integer $n\ge -1$ such that the modules $\H^i_\fa(M)$ are artinian
for all $i>n$ and all finite $R$-modules $M$. It is well-known
that all the local cohomology modules $\H^i_\fm(M)$ are artinian
and so $\q(\fm, R)=-1$.

Now recall the following well-known fact that ``a non-zero
$R$-module $M$ is artinian if and only if $\Supp M=\{\fm\}$ and
the socle of $M$, $\Hom_R(R/\fm,M)$, is finite''. Thus to
investigate artinianness of local cohomology module, one should
deal with the following two questions:

(1) When is a local cohomology module supported at the maximal
ideal?

(2) When is the socle of a local cohomology module finite?

The following conjecture, due to Huneke [{\bf Hu}], is related to
the question (2).

\noindent{\bf Conjecture.} Let $(R,\fm)$ be a regular local ring
and let $\fa$ be an ideal of $R$. Then the $R$-module $\Hom(R/\fm,
\H^n_\fa(R))$ is finite for all $n\ge 0$.

In [{\bf HS}], it is shown that the conjecture is true for regular
local rings in characteristic $p>0$. Also, Hartshorne gave an
example that the conjecture is not true for modules in place of
the ring, cf. [{\bf Ha2}].

For a finite $R$-module $M$ we define $\q(\fa, M)$ as the least
integer $n\ge -1$ such that the modules $\H^i_\fa(M)$ are
$\fm$-cofinite for all $i>n$. Note that the top local cohomology
$\H^{\dim M}_\fa(M)$ is artinian and so we have $\q(\fa, M)<\dim
M$. It is shown that for finite $R$-module $M$, the integer
$\q(\fa, M)$ is equal to the supremum of $\q(\fa, R/\fp)$ where
$\fp$ runs over the set $\Supp(M)$. One of the main result is as
follows:

Let $i$ be a given non-negative integer such that
$\H^i_\fa(R/\fb)$ is $\fm$-cofinite for all ideals $\fb$ of $R$.
Then for any $R$-module $M$, the integer $\q(\fa, M)$ is less than
$i$.

\vspace{.3in}

\newpage

\noindent{\bf 2. Finiteness of associated primes}

\vspace{.2in}

Grothendieck (compare [{\bf G}; Expos\'{e} XIII, 1.1]) proposed
the following conjecture.

``Let $M$ be a finite $R$-module and let $\fa$ be an ideal of $R$.
Then $\Hom_R(R/\fa, \H^j_{\fa}(M))$ is finite $R$-module for all
$j\ge 0$.''

\noindent Although this conjecture is not true in general, cf.
[{\bf Ha2}; Example 1], we have the following result.

\vspace{.1in}

\noindent{\bf Theorem 2.1} Let $\fa$ be an ideal of a noetherian
ring $R$. Let $s$ be a non-negative integer. Let $M$ be an
$R$-module such that $\Ext^i_R(R/\fa, M)$ is a finite $R$-module
for every $i\le s$, for example $M$ might be a finite $R$-module.
If $\H^i_\fa(M)$ is $\fa$-cofinite for all $i<s$, then
$\Hom_R(R/\fa, \H^s_\fa(M))$ is finite.

\vspace{.1in}

\noindent{\it Proof.} We use induction on $s$. If $s=0$, then
$\H^0_\fa(M)\cong\Gamma_\fa(M)$ and $\Hom_R(R/\fa,\Gamma_\fa(M))$
is equal to the finite $R$-module $\Hom_R(R/\fa, M)$. Now the
assertion holds.

Suppose that $s>0$ and the case $s-1$ is settled. Since
$\Gamma_\fa(M)$ is $\fa$-cofinite we have $\Ext^i_R(R/\fa,
\Gamma_\fa(M))$ is finite for all $i$. Now by using the exact
sequence $0\to \Gamma_\fa(M)\to M\to M/\Gamma_\fa(M)\to 0$ we get
$\Ext^i_R(R/\fa, M/\Gamma_\fa(M))$ is finite for all $i\le s$. On
the other hand $\H^0_\fa(M/\Gamma_\fa(M))=0$ and
$\H^i_\fa(M/\Gamma_\fa(M))\cong\H^i_\fa(M)$ for all $i>0$. Thus we
may assume that $\Gamma_\fa(M)=0$. Let $E$ be an injective hull of
$M$ and put $N=E/M$. Then $\Gamma_\fa(E)=0$ and $\Hom_R(R/\fa,
E)=0$. Consequently $\Ext^i_R(R/\fa, N)\cong\Ext^{i+1}_R(R/\fa,
M)$ and $\H^i_\fa(N)\cong\H^{i+1}_\fa(M)$ for all $i\ge 0$. Now
the induction hypothesis yields that $\Hom_R(R/\fa,
\H^{s-1}_\fa(N))$ is finite and hence $\Hom_R(R/\fa, \H^s_\fa(M))$
is finite too.\hfill$\square$

\vspace{.2in}

In proving finiteness of $\Hom_R(R/\fa, \H^1_\fa(M))$, we can even
give the following result with weaker condition.

\vspace{.1in}

\noindent{\bf Proposition 2.2} Assume $M$ is an $R$-module with
$\Ext^1_R(R/\fa, M)$ and $\Ext^2_R(R/\fa, \Gamma_\fa(M))$ are
finite are modules. Then $\Hom_R(R/\fa, \H^1_\fa(M))$ is finite.

\vspace{.1in}

\noindent{\it Proof.} By using the exact sequence
$0\to\Gamma_\fa(M)\to M\to M/\Gamma_\fa(M)\to 0$ we have the exact
sequence $\Ext^1_R(R/\fa, M)\to \Ext^1_R(R/\fa,
M/\Gamma_\fa(M))\to \Ext^2_R(R/\fa, \Gamma_\fa(M)).$ Thus the
module $ \Ext^1_R(R/\fa, M/\Gamma_\fa(M))$ is finite.

On the other hand the exact sequence $0\to M/\Gamma_\fa(M)\to
\mathrm{D}_\fa(M)\to \H^1_\fa(M)\to 0$, where
$\mathrm{D}_\fa(-)\cong\underrightarrow{\lim}\Hom_R(\fa^n, -)$,
induces the exact sequence
$$\Hom_R(R/\fa, \mathrm{D}_\fa(M))\to \Hom_R(R/\fa,
\H^1_\fa(M))\to\Ext^1_R(R/\fa, M/\Gamma_\fa(M)).$$ Since the left
side is zero and the right side is finite we get the
result.\hfill$\square$

\vspace{.2in}

The next result has been shown, using a spectral sequence
argument, by Divaani-Aazar and Mafi in [{\bf DiM}].

\vspace{.1in}

\noindent{\bf Corollary 2.3} Let $M$ be an $R$-module with
$\Ext^i_R(R/\fa, M)$ is finite for every $i$. Then the first non
$\fa$-cofinite local cohomology module of $M$ with respect to
$\fa$ has only finitely many associated primes.

\vspace{.2in}

The next result has been shown by Brodmann and Lashgari in [{\bf
BL}; Theorem 2.2] and it re-proves first part of  (2.4) and
generalizes (2.5)  of Brodmann, Rotthaus and Sharp [{\bf BRS} ].

\vspace{.1in}

\noindent{\bf Corollary 2.4} Let $M$ be a finite $R$-module. Let
$s$ be a non-negative integer such that $\H^i_\fa(M)$ is finite
for all $i<s$. Then the set $\Ass_R(\H^s_\fa(M))$ is finite.

\vspace{.3in}

\noindent{\bf 3. Cofiniteness}

Let $(R,\fm)$ be a local noetherian ring. Let $\fa$ be an ideal of
$R$. An $R$-module $M$ is $\fm$-cofinite if and only if
$\Supp(M)\subseteq\V(\fm)$ and $\Hom_R(R/\fm,M)$ is a
finite-dimensional vector space. As $\fm$-cofinite modules form an
abelian category, stable under taking submodules, quotient modules
and extensions, for an exact sequence $T_1\to T\to T_2$ of
$R$-modules, $T$ is $\fm$-cofinite provided $T_1$ and $T_2$ are
both $\fm$-cofinite.

\vspace{.2in}

\noindent{\bf Definition 3.1} Let $M$ be a finite $R$-module and
$\fa$ be an ideal of $R$. We define $\q(\fa, M)$ as the supremum
of the integers $i$ such that the module $\H^i_{\fa}(M)$ is not
$\fm$-cofinite.

\vspace{.1in}

As mentioned in the introduction $\q(\fa, M)<\dim(M)$. In the
following theorem it is shown that the invariant $\q(\fa, M)$
depends only on support of $M$.

\vspace{.1in}

\noindent{\bf Theorem 3.2} Let $\fa$ be a proper ideal of $R$ and
let $M$ and $N$ be finite $R$-modules such that $\Supp
N\subseteq\Supp M$. Then $\q(\fa, N)\le\q(\fa, M)$.

\vspace{.1in}

\noindent{\it Proof.} It is enough to show that $\H^i_\fa(N)$ is
$\fm$-cofinite for all $i>\q(\fa, M)$. We prove it by descending
induction on $i$ with $\q(\fa, M)<i\le\dim(M)+1$. For $i=\dim M+1$
there is nothing to prove. Now suppose $\q(\fa, M)<i\le\dim M$. By
Gruson's theorem [{\bf V}], there is a chain $$0=N_0\subseteq
N_1\subseteq\cdots\subseteq N_t=N$$ such that the factors
$N_j/N_{j-1}$ are homomorphic images of a direct sum of finitely
many copies of $M$. By using short exact sequences the situation
can be reduced to the case $t=1$. Therefore for some positive
integer $n$ and some finite $R$-module $L$ there exists an exact
sequence $$0\to L\to M^n\to N\to 0.$$ Thus we have the following
long exact sequence $$\cdots\to \H^i_{\fa}(L)\to\H^i_\fa(M^n)\to
\H^i_{\fa}(N)\to \H^{i+1}_{\fa}(L)\to\cdots.$$ By the induction
hypothesis $\H^{i+1}_{\fa}(L)$ is $\fm$-cofinite. Since
$\H^i_{\fa}(M^n)$ is $\fm$-cofinite, we have that $\H^i_{\fa}(N)$
is $\fm$-cofinite too.\hfill$\square$

\vspace{.2in}

The next Corollary gives a formula for $\q(\fa,-)$ in an exact
sequence.

\vspace{.1in}

\noindent{\bf Corollary 3.3} Let $0\to L\to M\to N\to 0$ be an
exact sequence of finite $R$-modules. Then $$\q(\fa,
M)=\max\{\q(\fa, L), \q(\fa, N)\}.$$

\vspace{.1in}

\noindent{\it Proof.} ``$\ge$'' holds by Theorem 3.2, and the
other side follows from the fact that $\fm$-cofinite modules form
an abelian category and from the long exact sequence in the proof
of Theorem 3.2.\hfill$\square$

\vspace{.2in}

\noindent{\bf Corollary 3.4} For an ideal $\fa$ of $R$ the
following holds
$$\q(\fa, R)=\sup\{\q(\fa, N)| N \mbox{is a finite
$R$-module}\}.$$

\vspace{.2in}

In the following theorem, we obtain that, for any finite
$R$-module, $\q(\fa, M)=\q(\fa, R/\fp)$ for some $\fp\in\Supp(M)$.

\vspace{.1in}

\noindent{\bf Theorem 3.5} Let $M$ be a finite $R$-module. Then
the following holds $$\q(\fa, M)=\sup\{\q(\fa, R/\fp)|\fp\in\Supp
M\}.$$

\vspace{.1in}

\noindent{\it Proof.} By Theorem 3.4 we have $\q(\fa,
R/\fp)\le\q(\fa, M)$ for all $\fp\in\Supp M$. Now assume that the
inequality is strict for all $\fp\in\Supp(M)$. There is a prime
filtration $0=M_0\subseteq M_1\subseteq\cdots\subseteq M_n=M$ of
submodules of $M$, such that for each $i$, $M_i/M_{i-1}\cong
R/\fp_i$ where $\fp_i\in\Supp M$. Set $t=\q(\fa, M)$. We have
$\H^t_\fa(R/\fp_i)$ is $\fm$-cofinite for all $1\le i\le n$. Thus
from the exact sequences $\H^t_\fa(M_{i-1})\to \H^t_\fa(M_i)\to
\H^t_\fa(R/\fp_i)$, $i=1,2,\cdots, n$, we eventually get $\q(\fa,
M_1)\ge t$ which is a contradiction.\hfill$\square$

\vspace{.2in}

\noindent{\bf Remark 3.6} In results 3.2--3.5 the underline ring
is assumed to be noetherian local. Now assume that $R$ is a
noetherian (not necessarily local) ring of finite dimension, $\fa$
is an ideal of $R$, and $M$ is a finite $R$--module. For each $\fp
\in {\Supp}M$, one can consider the invariant $\q(\fa R_{\fp},
M_{\fp})$ and define
$$ \q(\fa, M)={\sup}\{\q(\fa R_{\fp}, M_{\fp})|\fp \in {\Supp}M \}.$$
It is routine to check that all results 3.2--3.5 remain valid.

\noindent{\bf Remark 3.7} Let $\fa$ and $\fb$ be two ideals of
$R$. One can define $\q(\fa, \fb , M)$ as the supremum of the
integers $i$ such that the module $\H^i_{\fa}(M)$ is not
$\fb$-cofinite.

\noindent On the other hand we know that the category of
$\fb$-cofinite modules is an abelian subcategory of the category
of $R$-modules if one of the following holds:

\begin{verse}

(1) $R$ is complete local ring and $\fb$ is a one dimensional
prime ideal, cf. [{\bf DeM}; Theorem 2].

(2) $R$ is notherian ring of dimension at most two, cf. [{\bf M2};
Theorem 7.4].

\end{verse}
Therefore, with the same proof, we can see that the results
3.2--3.5 are valid for $\q(\fa, \fb, M)$ instead of $\q(\fa, M)$
provided (1) or (2) holds.

\vspace{.2in}

Now we are ready to give the main result of this section.

\vspace{.1in}

\noindent{\bf Theorem 3.8} Let $\fa$ be an ideal of $R$ and $i\ge
0$ be a given integer such that $\H^i_\fa(R/\fb)$ is
$\fm$-cofinite for all ideals $\fb$ of $R$. Then $\q(\fa,
R/\fp)<i$ for all $\fp\in\Spec R$. In particular, $\q(\fa, M)<i$
for all finite $R$-module $M$.

\vspace{.1in}

\noindent{\it Proof.} We prove by induction on $j\ge i+1$ that
$\H^j_\fa(R/\fp)$ is $\fm$-cofinite, for all $\fp\in\Spec R$. It
is enough to prove the case $j=i+1$. Suppose that
$\H^{i+1}_\fa(R/\fp)$ is not $\fm$-cofinite for some $\fp\in\Spec
R$, so that $\fa\nsubseteq\fp$. We first show that
$\Supp(\H^{i+1}_\fa(R/\fp))\subseteq\V(\fm)$. On the contrary,
assume $x\in\H^{i+1}_\fa(R/\fp)$ be a non-zero element whose
support is not in $\V(\fm)$. Since $x$ is annihilated by some
power of $\fa$, there exists $b\in\fa\setminus\fp$ such that
$bx=0$. Now consider the exact sequence $$0\lo R/\fp
\stackrel{b}{\lo} R/\fp\to R/({\fp+bR})\to 0$$ which induces the
following exact sequence
$$\H^i_\fa(R/({\fp+bR}))\to\H^{i+1}_\fa(R/\fp)\stackrel{b}{\lo}\H^{i+1}_\fa(R/\fp).$$

This shows that $\H^i_\fa(R/{\fp+bR})$ does not have support in
$\V(\fm)$, which, by the assumption, is impossible.

To show that $\H^{i+1}_\fa(R/\fp)$ is $\fm$-cofinite, we must show
in addition that $\Hom_R(R/\fm, \H^{i+1}_\fa(R/\fp))$ is finite-
dimensional vector space. Let $K$ be the kernel of multiplication
by $b$ in $\H^{i+1}_\fa(R/\fp)$, and consider the exact sequence
$$0\to K\to
\H^{i+1}_\fa(R/\fp)\stackrel{b}{\lo}\H^{i+1}_\fa(R/\fp).$$ Since
$b\in\fm$, it is easy to see that $\Hom_R(R/\fm,
K)=\Hom_R(R/\fm,\H^{i+1}_\fa(R/\fp))$. But $K$ is a quotient of $
\H^i_\fa(R/{\fp+bR})$, which is $\fm$-cofinite by hypothesis, so
$K$ is also $\fm$-cofinite. That means
$\Hom_R(R/\fm,\H^{i+1}_\fa(R/\fp))$ is of finite
dimension.\hfill$\square$

\vspace{.2in}

\noindent{\bf Remark 3.9} Let $R$ be a semi local ring with the
maximal ideals $\fm_1, \fm_2,\cdots ,\fm_n$. Let $J=\cap\fm_i$.
Then it is easy to check that the category of $J$-cofinite modules
is an abelian subcategory of the category of $R$-modules and the
results 3.2-3.6 are valid for $\q(\fa, J, M)$ instead of $\q(\fa,
M)$.

\vspace{.3in}

\noindent {\large\bf Acknowledgment.} The authors would like to
thank professor Divaani-Aazar for his useful comments on the last
section.

\baselineskip=16pt

\begin{center}
\large {\bf References}
\end{center}
\vspace{.2in}

\begin{verse}

[BH] W. Bruns; J. Herzog, {\em Cohen-Macaulay rings}, Cambridge
Studies in Advanced Mathematics, {\bf 39}. Cambridge University
Press, Cambridge, 1993.

[BL] M. P. Brodmann; A. Lashgari Faghani, {\em A finiteness result
for associated primes of local cohomology modules}, Proc. Amer.
Math. Soc., {\bf 128} (2000), 2851--2853.

[BRS] M. P. Brodmann; Ch. Rotthaus ; R. Y. Sharp, {\em On
annihilators and associated primes of local cohomology modules},
J. Pure Appl. Algebra {\bf 153} (2000), 197--227.

[DeM]D. Delfino; T. Marley, {\em Cofinite modules and local
cohomology}, J. Pure Appl. Algebra {\bf 121} (1997), no. 1,
45--52.

[DiM] K. Divaani-Aazar; A. Mafi, {\em Associated primes of local
cohomology modules}, Proc. Amer. Math. Soc., to appear.

[G] A. Grothendieck, {\em Cohomologie locale des faisceaux
coh\'{e}rents et th\'{e}or\'{e}mes de Lefschetz locaux et globaux}
$(SGA$ $2)$ Advanced Studies in Pure Mathematics, Vol. 2.
North-Holland Publishing Co., Amsterdam; Masson \& Cie, Editeur,
Paris, 1968.

[Ha1] R. Hartshorne, {\em Cohomological dimension of algebraic
varieties}, Ann. of Math. {\bf 88} (1968) 403--450.

[Ha2] R. Hartshorne, {\em Affine duality and cofiniteness},
Invent. Math. {\bf 9} (1970) 145--164.

[HL] C. Huneke; G. Lyubeznik, {\em On the vanishing of local
cohomology modules}, Invent. Math. {\bf 102} (1990), 73--93.

[HS] C. Huneke; R. Sharp, {\em Bass numbers of local cohomology
modules}, Trans. Amer. Math. Soc., {\bf 339} (1993), 765--779.

[Hu] C. Huneke, {\em Problems on local cohomology}, Free
resolutions in commutative algebra and algebraic geometry
(Sundance, UT, 1990), 93--108, Res. Notes Math., 2, Jones and
Bartlett, Boston, MA, 1992.

[KS] K. Khashyarmanesh; Sh. Salarian, {\em On the associated
primes of local cohomology modules}, Comm. Algebra {\bf 27}
(1999), 6191--6198.

[L] G. Lyubeznik, {\em Finiteness properties of local cohomology
modules (an application of D-modules to commutative algebra)},
Invent. Math. {\bf 113} (1993), 41--55.

[M1] L. Melkersson, {\em On asymptotic stability for sets of prime
ideals connected with the powers of an ideal}, Math. Proc. Camb.
Phil. Soc. {\bf 107} (1990), 267--271.

[M2] L. Melkersson, {\em Modules cofinite with respect to an
ideal}, Preprint, Link\"{o}pings University, 2003.

[M3] L. Melkersson, {\em Problems and results on cofiniteness: a
survey}, IPM Proceedings Series No. II, IPM 2004.

[S] A. K. Singh, {\em p-torsion elements in local cohomology
modules}, Math. Res. Lett. {\bf 7} (2000), 165--176.

[V] W. Vasconcelos, {\em Divisor theory in module categories},
North-Holland, Amsterdam, 1974.

\end{verse}

E-mail addresses:

\hspace{1.2in} dibaeimt@ipm.ir

\hspace{1.2in} yassemi@ipm.ir

\end{document}